\documentclass{amsart}

\usepackage{enumitem}
\usepackage[T1]{fontenc}
\usepackage[utf8]{inputenc}
\usepackage[english]{babel}
\usepackage{amsmath}
\usepackage{amsthm}
\usepackage{amssymb}
\usepackage[a4paper, margin=1in]{geometry}
\usepackage{comment}

\usepackage{subcaption}

\usepackage{graphicx}
\usepackage{tikz}
\usetikzlibrary{shapes.misc}
\usepackage{float}

\usepackage[unicode=true,pdfusetitle,
bookmarks=true,bookmarksnumbered=false,bookmarksopen=false,
breaklinks=false,pdfborder={0 0 1},backref=false,hidelinks]
{hyperref}

\makeatletter

\numberwithin{equation}{section}
\numberwithin{figure}{section}

\makeatother

\newtheorem{lemma}{Lemma}[section]

\begin{document}
	
	\title[Numerical methods for flow structure interactions]{Numerical methods for solving the linearized model of a hinged-free reduced plate arising in flow structure interactions
	}
	
	\author{
		Raj Narayan Dhara$^{1,2}$ 
	}
	
	\author{
		Krzysztof E. Rutkowski$^{2}$ 
	}

	\author{
		Katarzyna Szulc$^{2}$ 
	}
	\thanks{$^1$ \textit{KMA, Faculty of Science}, \textit{Palack\'y University Olomouc}, Czech Republic.
	}
	\thanks{$^2$ \textit{Systems Research Institute} \textit{Polish Academy of Sciences}, Warsaw, Poland
	}
	
	\maketitle
	
	\begin{abstract}
		The problem of partially hinged and partially free rectangular plate that aims to represent a suspension bridge subject to some external forces (for example the wind) is considered in order to model and simulate the unstable end behavior.
		Such a problem can be modeled by a plate evolution equation, which is nonlinear with a nonlocal stretching effect in the spanwise direction.
		The external forces are periodic in time and cause the vortex shedding on the structure (on the surface of the plate) and thus it may cause damage to the material.
		Numerical study of the behavior of steady state solutions for different values of the force velocity are provided with two finite element methods of different type.
	\end{abstract}
	\section{Introduction}
	The general problem that we are interested in is the partially hinged and partially free rectangular plate that aims to represent a suspension bridge subject to some external forces (for example the wind).
	We want to model and simulate the unstable end behavior.
	Such a problem can be modeled by a plate evolution equation, which is nonlinear with a non-local stretching effect in the span-vise direction.
	Originally, the flow of the wind in the chord-wise direction was modeled through a so-called piston-theoretic approximation. Such an approximation provides both weak dissipation and non-conservative forces.
	We suppose the external forces are periodic in time and cause the vortex shedding on the structure (on the surface of the plate). 
	The latter leads to a damage to the structure.
	Therefore, we want to see what happens in the stationary case because under some conditions, the solution to the boundary value problem under consideration can be controllable to a single stationary point, which is trivial.
	
	The main objective of this paper is to study numerically the behavior of steady state solutions for different values of the force velocity.
	These are related to asymptotic, long time behaviour of evolutions driven by \eqref{eq:2.1}.
	The force velocity given by the parameter $\alpha$ is a key player from the point of view of the non-self-adjoint stationary problem \eqref{eq:8.1}.
	Indeed, if it is small, the uniqueness of a stationary solution is guaranteed. Otherwise, the existence of multiple unimodal stationary solutions is asserted by the Theorem 3.5 of \cite{lasweb} and their number grows with the parameter $\alpha$.
	This means that we deal with multiple sets of stationary solutions, which provide different behavior for long time dynamics and, in consequence, emphasize the complexity of the global attractor.

	The paper is organized in the following way.
	The model of a hinged-free nonlinear and nonlocal boundary value problem for the thin rectangular plate is defined, and the variational formulation is recalled in Section \ref{chapter:model}.
	In Section \ref{chapter:stationary_case}, the stationary case is considered.
	Section \ref{chapter:numerical_rectangular} and \ref{chapter:numerical} is entirely dedicated to the description of the finite element approximation and the detailed algorithm for the numerical method is described.
	Some estimates are provided for convergence, and finally, numerical results,  in Section \ref{chapter:numerical}, are presented.
	The last chapter concludes the present work.

	\section{The model}\label{chapter:model}
	
	The original model for flow structure interaction describes the interaction between a nonlinear plate with the flow above it. Hence the plate should be treated two-dimensionally so the equation describing the downwash of the flow and the oscillations of the plate can agree on the interface. 
	The behavior of the gas is modeled via the theory of potential flows \cite{dowell} which produces a convected wave equation for the perturbed velocity potential of the flow.
	The oscillatory behavior of the plate is governed by the Kirchhoff plate equation with the von Karman nonlinearity, which is used in the modeling of the large oscillations of thin, flexible plates.
	This goes back to the fully coupled flow-plate model considered in~\cite[eq. (2.6)]{CLW2014} which was reduced to a delayed von Karman plate~\cite[eq. (2.7)]{CLW2014}. 
	
	In the model~\eqref{eq:2.1} that we present in this paper, 
	the aerodynamical pressure is considered as $p({\bf x},t)=g({\bf x}) + (\partial_t +\alpha \partial_y){\rm tr}[\phi]$, where ${\rm tr}[\phi]=\phi|_{z=0}$ is the trace of the flow $\phi$ on the plane $\{{\bf x}\in\mathbb{R}^3 : z= 0\}$. 
	In contrast to the original model, {\em no delay} is taken here.
	The fluid flow environment we consider is $\mathbb{R}^{3}_{+}=\{(x,y,z):\ z>0 \}$. 
	The plate is immersed in an inviscid flow (over body, $z\ge 0$) with velocity $\alpha\neq 1$ in the $y$-direction. 
	The damping in this model is given by the linear, frictional damping $(k+1)u_t$ with $k \ge -1$. 
	Furthermore, the strong coupling here takes place in the downwash term of the flow potential (the Neumann condition in~\cite[eq. (2.4)]{CLW2014}) by $d({\bf x},t)=(\partial_t +\alpha \partial_y)u(x,y)\cdot 1_{\Omega}({\bf x}
	),\ (x,y)\in\mathbb{R}^2$. 
	In particular, we write~\eqref{eq:2.1} as follows
	\begin{align}\label{eq:2:1a}
		u_{tt}+(k+1)u_t +\Delta^2 u+\left[P-S\int_{\Omega}u_x^2\, dx\right]u_{xx}=g+(\partial_t +\alpha \partial_y)u \quad \mathrm{in} \ \Omega\times(0,T).
	\end{align}
	
	\begin{figure}[htp]
		\begin{center}
			\includegraphics[width=0.3\textwidth]{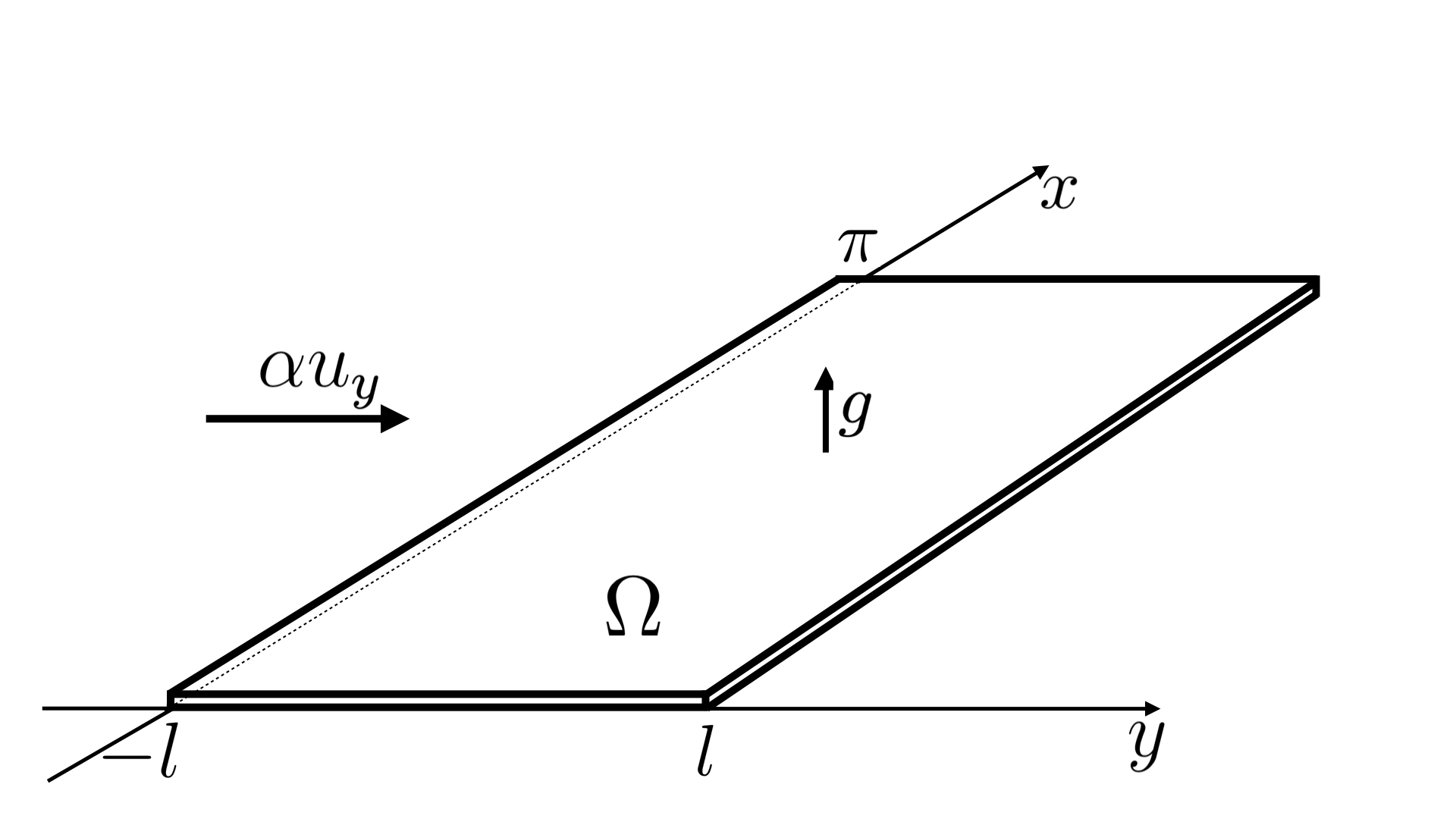}
			\caption{Rectangular plate representing the domain $\Omega$.}\label{fig1}
		\end{center}
	\end{figure}

		We focus here on a rectangular plate as presented on Fig.\ref{fig1}.
		The plate occupies the domain $\Omega=(0,\pi)\times(-l,l)$.
		We suppose that the edges $\{0\}\times[-l,l]$ and $\{\pi\}\times[-l,l]$ are hinged and the edges $[0,\pi]\times\{-l\}$ and $[0,\pi]\times\{l\}$ are free.
		We denote by $\Gamma_D\cup\Gamma_N$ the boundary of $\Omega$ with $\Gamma_D=\{0,\pi\}\times[-l,l]$ and $\Gamma_N=[0,\pi]\times\{-l,l\}$.
		The movement of the hinged-free plate is modeled by the nonlinear and nonlocal evolution equation.
		Such an equation was introduced in \cite{12,21}.
		In this paper, we adapt the description of this model from \cite{lasweb} and also recall the necessary notations for the equations under consideration.
		We start with the strong formulation of the displacement $u=u(x,y,t)$ given in the following way:
		\begin{equation}\label{eq:2.1}
			\left\{
			\begin{array}{rcl}
				u_{tt}+ku_t+\Delta^2 u+[P-S\int_{\Omega}u_x^2]u_{xx}
				=g+\alpha u_y  & \mathrm{ in }&  \Omega\times(0,T)\\
				u=u_{xx}=0  &\mathrm{ on }&  \Gamma_D\\
				u_{yy}+\sigma u_{xx}= u_{yyy}+(2-\sigma)u_{xxy}=0 & \mathrm{ on }&   \Gamma_N\\
				u(x,y,0)=u_0(x,y),\,\, u_t(x,y,0)=v_0(x,y)  &\mathrm{ in }&  \Omega.
			\end{array}
			\right.
		\end{equation}
		Here, $k$ is an overall damping, $\alpha\in\mathbb{R}$ is a generalized flow parameter corresponding to the force of the wind, and $\sigma\in(0,1)$ is a Poisson ratio.
		Boundary conditions correspond to the longitudinally hinged, laterally free plate.
		Parameter $P$ is the so-called prestressing parameter.
		
		In this paper, especially in the numerical computations, we mostly deal with the case of a strongly prestressed plate, i.e., when the parameter $P$ comes from the range $0<P<\lambda_1$ where $\lambda_1$ is a first eigenvalue of the principal structural operator defined below in \eqref{eq:2.3}.
		However, it is possible to also consider the case of a weakly prestressed plate when considering stability. In this case, $P$ is taken from the range $\lambda_1\leq P<\lambda_2$, with $\lambda_2$ the second eigenvalue of the principal structural operator.
		However, most of the result from \cite{lasweb} allows $P\in \mathbb{R}$.
		The parameter $S\geq0$ indicates the strength of the restoring forces resulting from the stretching in $x$.
		
		Let us denote by $H^s(\Omega)$ the Hilbert Sobolev space of order $s\in\mathbb{N}$ with norm $\|\cdot \|_s$ and semi norm given by
		\begin{equation}\label{norma1} 
			\| u \|_s=\left(\sum_{|\nu|\leq s}\int_{\Omega}|D^{s}u|^{2}\right)^{1/2}, \quad
			| u |_s=\left(\sum_{|\nu|= s}\int_{\Omega}|D^{s}u|^{2}\right)^{1/2}
		\end{equation}
		for all $u\in H^{s}(\Omega)$ with $\nu$ a multiindex such that $\nu=(\nu_1,\nu_2)$, $\nu\geq0$, $|\nu|=\nu_1+\nu_2$ and $\partial^{\nu}=(\frac{\partial}{\partial x})^{\nu_1}(\frac{\partial}{\partial y})^{\nu_2}$.
		The inner product in $L^2(\Omega)$ we write as $(\cdot,\cdot)$.
		We denote by $\mathcal{K}^2$ the phase space of admissible displacements $u$ for the partially clamped partially free plate from the equation \eqref{eq:2.1}. This space is given by
		\begin{equation}\label{eq:2.0}
			\mathcal{K}^2=\left\{ u\in H^2(\Omega):\, u=0\  \mathrm{ on }\ \Gamma_D \right\}
		\end{equation}
		with the following scalar product:
		\begin{equation}\label{eq:2.2}
			a(u,v):=\int_{\Omega}(\Delta u\Delta v-(1-\sigma)[u_{xx}v_{yy}+u_{yy}v_{xx}-2u_{xy}v_{xy}]),
		\end{equation}
		with $u,v\in \mathcal{K}^2(\Omega)$.
		Notice that the latter induces a norm 
		\begin{equation}\label{norm:K2}
			\|u\|_{\mathcal{K}^{2}}=\sqrt{a(u,u)}
		\end{equation}
		equivalent to the usual Sobolev norm $\|\cdot\|_{2}$.
		Next, we define the positive, self-adjoint biharmonic operator corresponding to $a(\cdot,\cdot$), $A:L^2(\Omega)\to L^2(\Omega)$.
		This operator, taken with the boundary conditions in \ref{eq:2.1} is given by $Au=\Delta^2 u$. Its domain we denote by
		\begin{align}\label{eq:2.3}
			\mathcal{D}(A)=\left\{ u\in H^4(\Omega)\cap \mathcal{K}^2 
			: u_{xx}=0 \ \mathrm{on}\   \Gamma_D,
			u_{yy}+\sigma u_{xx}=0, u_{yyy}+(2-\sigma)u_{xxy}=0  \ \mathrm{on}\    \Gamma_N
			\right\}
		\end{align}
		The conditions $u_{xx}=0$ on $\Gamma_D$ and $u_{yy}+\sigma u_{xx}=0$, $u_{yyy}+(2-\sigma) u_{xxy}=0$
		on $\Gamma_N$ are the natural conditions associated with $a(\cdot,\cdot)$ in the strong form.
		The characterization in (\ref{eq:2.3}) seems to not be obvious due to the presence ofht corners and the mixed boundary conditions occuring on the same boundary. 
		However, in the case of the plate equation defined in the rectangular domain where the angles of the corners are equal to $90^\mathrm{o}=\pi/2$ the absence of any singular behaviors in the bending moments at the corners can be verified, as in \cite{MANSFIELD}, Section 3.8.
		
		{\bf Weak formulation.}
		Suppose that $g\in L^2(\Omega)$.
		A weak solution of \eqref{eq:2.1} is a function 
		$$
		u\in C^0(\mathbb{R}_{+},\mathcal{K}^{2}(\Omega))\cap C^1(\mathbb{R}_{+},L^{2}(\Omega))\cap C^2(\mathbb{R}_{+},(\mathcal{K}^{2}(\Omega))^\prime)
		$$
		such that for all $t>0$ and all $v\in \mathcal{K}^{2}(\Omega)$, we have
		\begin{align}\label{eq:2.4}
			\langle u_{tt},v\rangle+k(u_t,v)+a(u,v)+[S\|u_x \|^2_0-P](u_x,v_x)
			=(g,v)+\alpha (u_y,v).
		\end{align}
		Note that strong solutions are generalized, and generalized solutions are weak \cite{n16}.
		
		\section{Stationary case and Linearisation}\label{chapter:stationary_case}

		For the numerical computations, we focus on the stationary case.
		As it was mentioned previously, in some cases the attractor may reduce to the unique stationary point, which can be considered as a {\em trivial solution}.
		In order to prove that the attractor, in particular, cannot be reduced to a single stationary solution (recall that $g\in L^2(\Omega)$),
		we seek solutions to the stationary case.
		Thus, we consider the following stationary case of the equation \eqref{eq:2.1}:
		\begin{equation}\label{eq:8.1}
			\left\{
			\begin{array}{rcl}
				\Delta^2 u+[P-S\int_{\Omega}u_x^2]u_{xx}=g+\alpha u_y & \mathrm{in} & \Omega\\
				u=0,\quad u_{xx}=0 & \mathrm{on} & \Gamma_D\\
				u_{yy}+\sigma u_{xx}= u_{yyy}+(2-\sigma)u_{xxy}=0 & \mathrm{on} &  \Gamma_N
			\end{array}
			\right.
		\end{equation}

		Notice that the bilaplacian operator given in \eqref{eq:2.2}
		induces a norm $\|u\|_{\mathcal{K}^{2}}=\sqrt{a(u,u)}$ which is equivalent to the usual Sobolev norm $\|\cdot\|_2$.
		One can show that for the stationary nonlinear problem \eqref{eq:8.1}, the solutions are bounded in the sense of its norm, i.e., they are taken from a ball of finite radius.
		This result is presented in the following lemma:
		\begin{lemma}\label{lem41}
			Let $S>0$ and $0\leq P<\lambda_1$ be given and suppose that $u$ is the solution to the problem \eqref{eq:8.1}.
			There exists constant $C>0$ depending on $\alpha$, $g$, $P$ and $S$ such that $u\in B_C(\mathcal{K}^2)$ where $B_C$ is a ball of radius $C$.
		\end{lemma}

	In this paper, we are interested, however, in nontrivial solution to the problem \eqref{eq:8.1} thus we introduce the related linearized problem:
	\begin{equation}\label{eq:8.2}
		\left\{
		\begin{array}{rcl}
			\Delta^2 U-\mu U_{xx}=G+\alpha U_y & \mathrm{in} & \Omega\\
			U=0,\quad U_{xx}=0 & \mathrm{on} &\Gamma_D\\
			U_{yy}+\sigma U_{xx}=U_{yyy}+(2-\sigma)U_{xxy}=0 & \mathrm{on} &  \Gamma_N
		\end{array}
		\right.
	\end{equation}
	where $\mu$ corresponds to the term $\left[S\displaystyle\int_{\Omega}u_x^2-P\right]$, $\alpha$ is  a flow parameter that represents the wind speed, $G$ is a positive constant and $\sigma\geq0$ measures the weak frictional damping. Let us formulate the following lemma
	\begin{lemma}\label{lem22}
		Suppose that $U\neq 0$ is the solution to the problem \eqref{eq:8.2} for $P\in\mathbb{R}$, $\mu>-P$ and $G=\frac{\|U_x\|_0}{\sqrt{\frac{\mu +P}{S}}}g$, $g$ is a function as in \eqref{eq:8.1}. 
		Then the function 
		\begin{equation}\label{eq:11}
			u(x,y)=\sqrt{\frac{\mu +P}{S}}\frac{U(x,y)}{\|U_x\|_0}
		\end{equation} 
		is a nontrivial solution to \eqref{eq:8.1}.
	\end{lemma}

	Now we define the weak solution for stationary problem \eqref{eq:8.1}.
	For all $v\in\mathcal{K}^{2}(\Omega)$, 
	the weak solution of \eqref{eq:8.1} is a function
	$$
	u\in C^0(\mathbb{R}_{+},\mathcal{K}^{2}(\Omega))\cap C^1(\mathbb{R}_{+},L^{2}(\Omega))
	$$
	such that 
	\begin{equation*}
		a(u,v)+[S\|u_x \|^2_0-P](u_x,v_x)=(g,v)+\alpha (u_y,v)
	\end{equation*}
	Based on this formula,
	we construct the following weak formulation for the linearized problem \eqref{eq:8.2} in the form:
	\begin{equation*}
		a(U,V)+\mu(U_x,V_x)=(G,V)+\alpha (U_y,V), \quad U,V\in \mathcal{K}^{2}(\Omega)
	\end{equation*}
	where $a(U,V)$ is given by:
	\begin{align}\label{eq14}
		a(U,V)= 
		\int_{\Omega}(\Delta U\Delta V+(1-\sigma)[2U_{xy}V_{xy}-U_{xx}V_{yy}-U_{yy}V_{xx}]),
	\end{align}
	with $U,V\in \mathcal{K}^2(\Omega)$.
	Finally, 
	we can formulate the variational form in the following way:
	\begin{equation}\label{vf4}
		\left\{
		\begin{array}{l}
			\text{Find } U\in \mathcal{K}^2(\Omega) \text{ such that }\\
			a(U,V) +\mu(U_x,V_x)=(G,V)+\alpha (U_y,V) \\ 
			\text{for all} \,V\in  \mathcal{K}^2(\Omega)
		\end{array}
		\right.
	\end{equation}
	with the domain of the operator $Au=\Delta^2 u$ given by \eqref{eq:2.3}.
	
	\section{Numerical analysis}
	In this section, we study the approximate solution to the linearized plate equation in a form \eqref{eq:8.2} obtained with the finite element method.
	We consider two cases for this analysis. 
	The first case bases on the model of solution in the form of separable variables (as it was suggested in \cite{gaz23}) and described on a rectangular grid.
	The second case bases on standard approach of P4 finite element method and is considered in  triangular discretization. The method is related to the approach presented in \cite{szulc} but in case of hinged boundary conditions on a part $\Gamma_d$ of the boundary.
	In both cases the geometry of the domain is a rectangular two dimensional plate given by $(0,\pi)\times(-0.2,0.2)$ with the Poisson ratio $\sigma=0.2$ and $\mu=-0.5$.
	
	\subsection{Rectangular grid case}\label{chapter:numerical_rectangular}

	We take into account that the solution should satisfy $u=u_{xx}=0$ on $\Gamma_D$.
	
	We denote:
	$M_1$ - number of mesh elements in $x$ direction,
	$M_2$ - number of mesh elements in $y$ direction,
	$N_1=M_1$ - number of mesh nodes in $x$ direction,
	$\overline{N}_1=N_1-2$ - number of interior mesh nodes in $x$ direction (i.e. $x_2,\dots,x_{N_1-1}$), 
	$N_2=3 M_2+1$ - number of mesh nodes in $y$ direction,
	$N=M_1\cdot (3M_2+1)$ - total number of mesh nodes,
	$M_1\cdot M_2 $ - total number of mesh elements,
	$\tilde{N}=\overline{N}\cdot N_2$ -  number of nodes that are not on the boundary $\Gamma_D$,
	$X_{s}:=(x_{i},y_{j})\ i = 1, \dots, N_1;\ j=1,\dots, N_2$ finite element nodes.
	First, we form a uniform rectangular partition of $\Omega:=(0,\pi)\times (-\ell,\ell)$ into $N_1$
	elements in $x$-axis and $N_2$ elements in $y$-axis.

	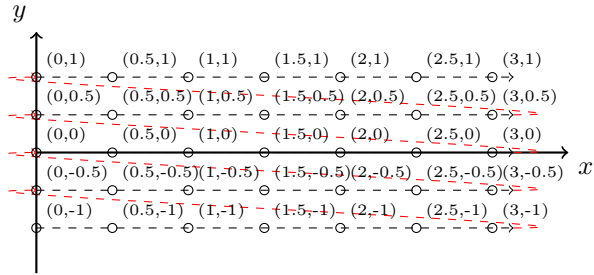
\begin{figure}[h]
		\centering
		\begin{tikzpicture}[scale=2]
			\draw[thick,->] (0,0) -- (3.5,0) node[anchor=north west] {$x$};
			\draw[thick,->] (0,-0.8) -- (0,0.8) node[anchor=south east] {$y$};
			\foreach \y in {-1,-0.5,...,1}
			\draw[dashed,->] (0,\y/2) -- (3.14,\y/2);
			\foreach \y in {-1,-0.5,...,0.5}
			\draw[red,dashed,->] (3.14,\y/2) to[in=180,out=0] (0,\y/2+0.25);
			\foreach \x in {0,0.5,...,3.14}
			{\foreach \y in {-1,-0.5,...,1}
				{
					\draw (\x,\y/2) circle (0.03);
					\node[above right] at (\x,\y/2) {\tiny (\x,\y)};
				}
			}
		\end{tikzpicture}
		\caption{Uniform Rectangular Mesh on $[0, \pi]\times [-1/2,1/2]$;}
		\label{fig:mesh2}
	\end{figure}
	The Galerkin formulation (without considering the Dirichlet
	boundary condition, which will be handled later)is formulated as follows : find $u_h\in U_h$ such
	that
	\begin{align}\label{eq:varf2}
		A(u_h, v_h) = (g , v_h),\ \text{for any}\  v_h\in U_h,
	\end{align}
	where
	\begin{equation*}
		A(u,v)=a(u,v)+\mu(u_x,v_x)-\alpha(u_y,v),\, \text{for all } u,v \in  \mathcal{K}^2(\Omega) 
	\end{equation*}
	and $U_h(x,y):= {\rm span}\{\Psi_i(x)\Phi_j(y) \}_{i=1,\dots, N_1, j=1,\dots, N_2}$ is chosen to be a finite element space where $$\{\Psi_i(x)\Phi_j(y) \}_{i=1,\dots, N_1, j=1,\dots,N_2}$$ are the global finite element basis functions to be defined as below.\\
	{\bf Construction of local basis functions:}
	Recall $N$ denote the number of local finite element nodes (local
	finite element basis functions) in a mesh element. Nodes of the mesh are in positions $(x_i,y_j)$ for $i=1,\dots,N_1$, $j=1,\dots,N_2$.

	Consider the following (shape) functions separately in $x$ and $y$-direction. 
	
	For $i=1,\dots \overline{N}_{1}$, let $\Phi_i(x)$ be Fourier shape functions with only non-zero coefficients related to the sinus functions to interior vertices $x_i$, $i=2,\dots,N_1-1$ (see \cite{highorder}):
	\begin{align*}
		\tilde{\Psi}_i(x)=\sin(i\cdot  x), \quad x \in [0,\pi], \ i=1,\dots, \overline{N}_1
	\end{align*}
	Note that the functions $\tilde{\Psi}_i$, $i=1,\dots, \overline{N}_1$, called basis of interpolants, does not have the Kronecker-delta-property at distributed points $x_i$, $i=2,\dots, N_1-1$, therefore one must use so-called transfer matrix $\mathcal{T}$ defined as
	\begin{equation*}
		\mathcal{T}_{i_1,i_2}=\Psi_{i_2}(x_{i_1+1}),\quad  i_1,i_2=1,\dots,\overline{N}_1
	\end{equation*}
	and define normalized basis as
	\begin{equation*}
		\Psi_i(x)=  \sum_{k=1}^{\overline{N}_1}  \tilde{\Psi}_i(x_{k+1}) \cdot (\mathcal{T}^{-1})_{k,i}.
	\end{equation*}

		\begin{figure}[H]
			\centering
			\begin{subfigure}{0.45\textwidth}
				\includegraphics[width=1\textwidth]{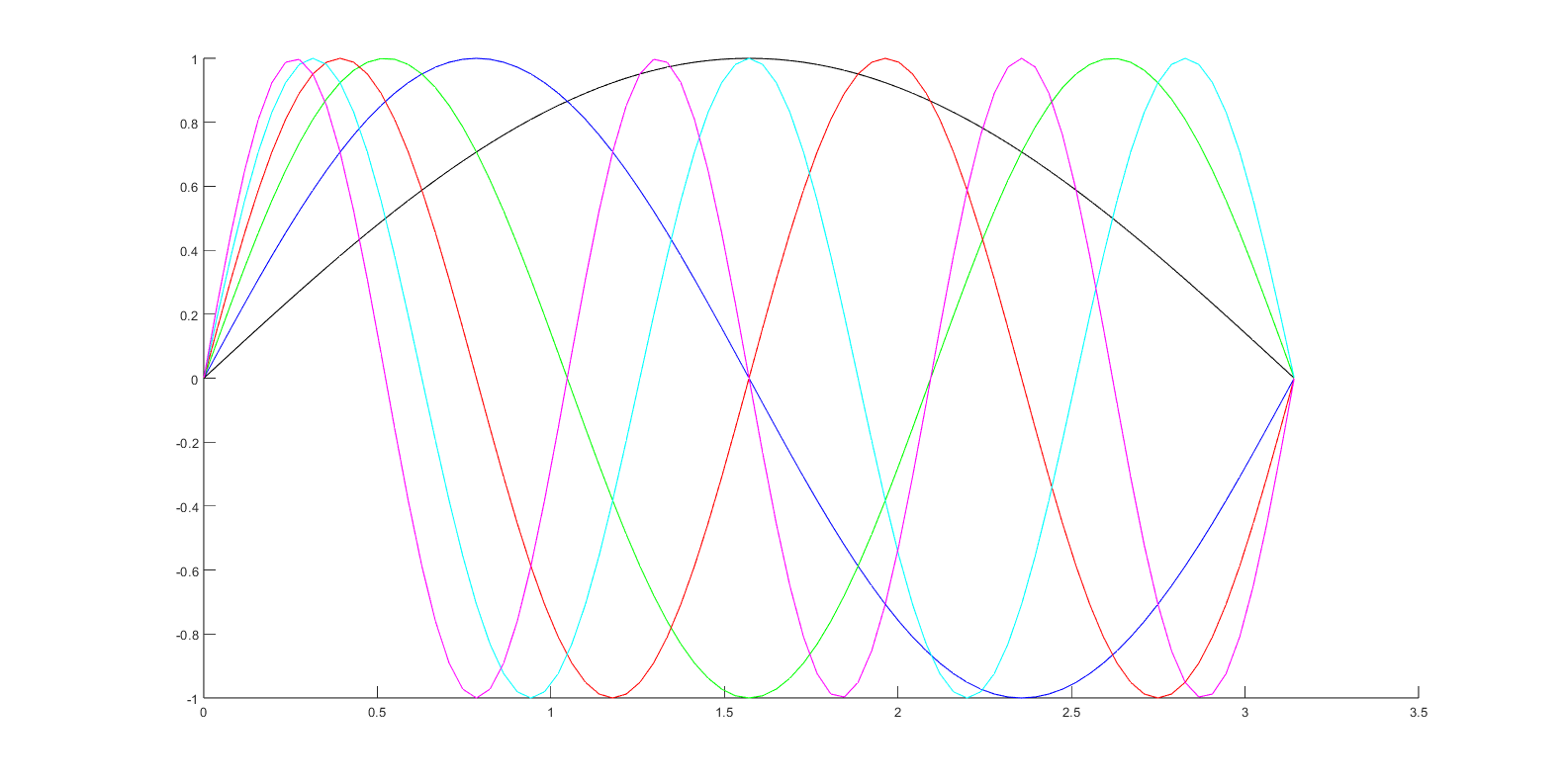}
				\caption{Unnormalized basis}\label{sin6n}
			\end{subfigure}
			\begin{subfigure}{0.45\textwidth}
				\includegraphics[width=1\textwidth]{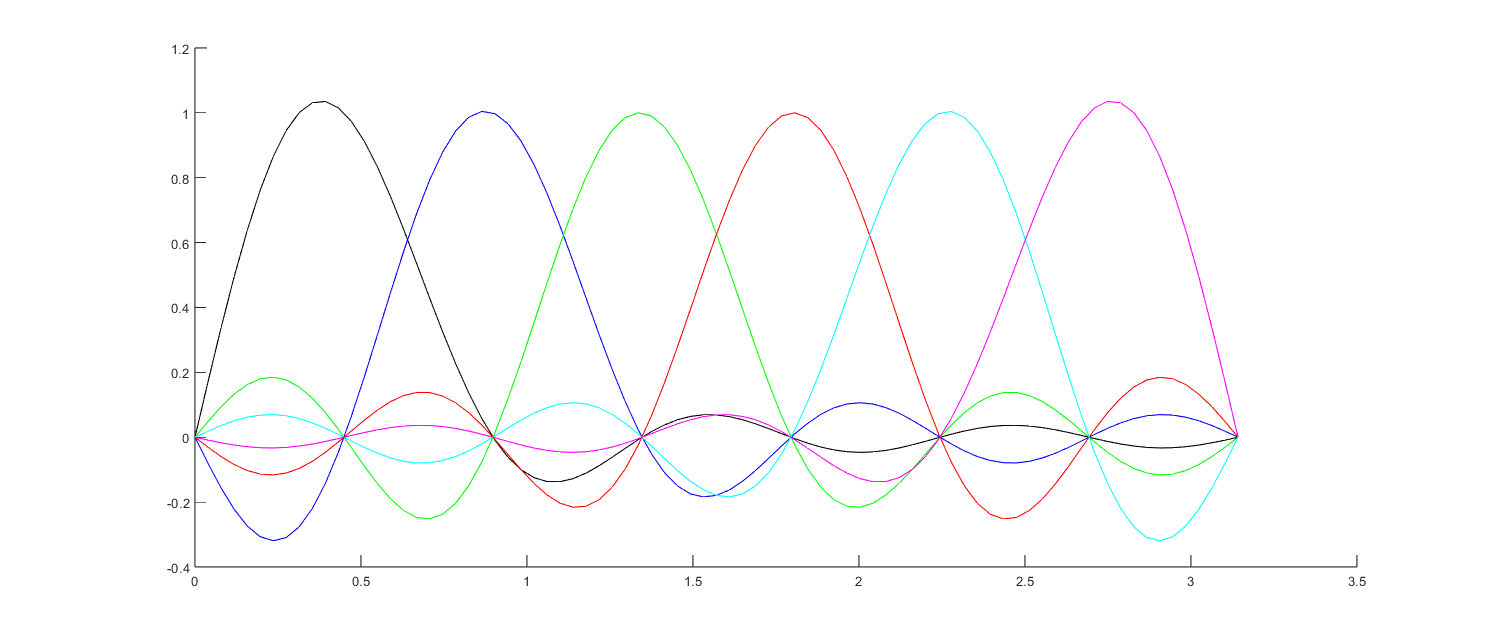}
				\caption{Normalized basis}\label{sin6}
			\end{subfigure}
			\caption{Basis functions for $x$ variable corresponding to 6 interior nodes}
			\label{fig:figures}
		\end{figure}

		The hinged boundary conditions are verified on $\Gamma_D$ (i.e. on boundary vertices $x_1$ and $x_{N_1}$) since the sinus function and its second derivative is equal to zero at $x=x_1=0$ and $x=x_{N_1}=\pi$.
		
		For $j=1,\dots,N_2$ let $\Psi_i(y)$ be Lagrangian functions of order 3 defined as:
		\begin{enumerate}
			\item for $j=3k+2,\ k=0,\dots,N_2-1$:
			\begin{align*}
				\phi_j(y):=&
				\left\{\begin{array}{cl}
					\frac{y-y_{j-1}}{y_{j}-y_{j-1}}\cdot \frac{y-y_{j+1}}{y_{j}-y_{j+1}} \cdot \frac{y-y_{j+2}}{y_{j}-y_{j+2}} 
					& y\in[y_{j-1},y_{j+2}], \\
					0, & \text{otherwise}, 
				\end{array}
				\right.
			\end{align*}
			\item for $j=3k+3,\ k=0,\dots,N_2-1 $:
			\begin{align*}
				\phi_j(y):=&
				\left\{\begin{array}{cl}
					\frac{y-y_{j-2}}{y_{j}-y_{j-2}}\cdot \frac{y-y_{j-1}}{y_{j}-y_{j-1}} \cdot \frac{y-y_{j+1}}{y_{j}-y_{j+1}}
					& y\in[y_{j-2},y_{j+1}], \\
					0, & \text{otherwise}, 
				\end{array}
				\right.
			\end{align*}
			\item for $j=3k+4,\ k=0,\dots,N_2-2$:
			\begin{align*}     
				\phi_j(y):=&
				\left\{\begin{array}{cl}
					\frac{y-y_{j-3}}{y_{j}-y_{j-3}}\cdot 
					\frac{y-y_{j-2}}{y_{j}-y_{j-2}}\cdot 
					\frac{y-y_{j-1}}{y_{j}-y_{j-1}}
					& y\in[y_{j-3},y_{j})\\
					\frac{y-y_{j+1}}{y_{j}-y_{j+1}}\cdot 
					\frac{y-y_{j+2}}{y_{j}-y_{j+1}}\cdot 
					\frac{y-y_{j+3}}{y_{j}-y_{j+3}} 
					& y\in[y_{j},y_{j+3}], \\
					0, & \text{otherwise}, 
				\end{array}
				\right. 
			\end{align*}
			\item for $j=1$
			\begin{align*}         
				\phi_1(y):=&
				\left\{\begin{array}{cl}
					\frac{y-y_{2}}{y_{1}-y_{2}}\cdot 
					\frac{y-y_{3}}{y_{1}-y_{3}}\cdot 
					\frac{y-y_{4}}{y_{1}-y_{4}} 
					& y\in[y_{1},y_{4}], \\
					0, & \text{otherwise}, 
				\end{array} \right.
			\end{align*}
			\item for $j=3\cdot N_2+1$
			\begin{align*} 
				\phi_{3\cdot N_2+1}(y):=&
				\left\{\begin{array}{cl}
					\frac{y-y_{3\cdot N_2-2}}{y_{3\cdot N_2+1}-y_{3\cdot N_2-2}}&
					\cdot \frac{y-y_{3\cdot N_2-1}}{y_{3\cdot N_2+1}-y_{3\cdot N_2-1}}\cdot \\
					\cdot \frac{y-y_{3\cdot N_2}}{y_{3\cdot N_2+1}-y_{3\cdot N_2}} ,
					& y\in[y_{3\cdot N_2-2},y_{3\cdot N_2+1}], \\
					0, & \text{otherwise}. 
				\end{array}
				\right.
			\end{align*}
		\end{enumerate}
		Since $u_h\in U_h= {\rm span}\{\Psi_i \Phi_j \}_{i=1,\dots, N_1, j=1,\dots,N_2}$ , then
		\begin{align*}
			u_h (x,y)= \sum_{i=1}^{N_1}\sum_{j=1}^{N_2} q_{ij}\Psi_i(x) \Phi_j(y)
		\end{align*}
		for some coefficients $q_{ij},\ i=1,\dots, N_1, j=1,\dots,3\cdot N_2$. We have that
		\begin{align}\label{basis:def2}
			\Psi_i\Phi_j\circ (X_t) =\delta_{tij}= \begin{cases}
				0,\ \text{if}\ t=(i_1,j_1)\neq (i,j)\\
				1,\ \text{if}\ t=(i_1,j_1)=(i,j);
			\end{cases}
		\end{align}                 
		hence with $t=(j-1)\cdot (N_1-1)+i$, $i=1,\dots,N_1,\ j=1,\dots, N_2$
		\begin{align}\label{basis:repr2}
			u_h(X_k) =& \sum_{t=1}^{\tilde{N}} q_t\Psi_i\Phi_j\circ (X_t) = q_k.
		\end{align}
		Therefore the coefficient $q_t$ is actually the numerical solution at the node $X_t\ (t=1,\cdots \tilde{N})$.  We choose the test function $v_h(x,y) = \Psi_{\tilde{i}}(x)\Phi_{\tilde{j}}(y)$, $\tilde{i}=1,\dots N_1$ $ \tilde{j}=1,\dots N_2$. Then the finite element formulation~\eqref{eq:varf2} gives: for $\tilde{i}=1,\dots N_1,\ \tilde{j}=1,\dots,N_2$ and $t=(j-1)\cdot (N_1-1)+i$
		\begin{equation}\label{eq:varf:2}
			A\left(\sum_{t=1}^{\tilde{N}} q_t\Psi_i \Phi_j, \Psi_{\tilde{i}}\Phi_{\tilde{j}}\right) = (g , \Psi_{\tilde{i}} \Phi_{\tilde{j}}),
		\end{equation}
		and for $\tilde{i}=1,\dots N_1,\ \tilde{j}=1,\dots, N_2$ $t=(j-1)\cdot (N_1-1)+i$
		\begin{equation}\label{eq:varf:21}
			\sum_{t=1}^{\tilde{N}}q_t A(\Psi_i \Phi_j, \Psi_{\tilde{i}} \Phi_{\tilde{j}} ) = (g , \Psi_{\tilde{i}} \Phi_{\tilde{j}} )
		\end{equation}
		{\bf Matrix formulation:} Define the stiffness matrix and the load vector,
		respectively,
		\begin{align*}
			\mathcal{A}&= [\mathcal{A}_{st}]_{s,t=1}^{\tilde{N}}:= A(\Psi_i \Phi_j, \Psi_{\tilde{i}} \Phi_{\tilde{j}})\quad {\bf b}:= [b]_{s=1}^{\tilde{N}} = (g , \Psi_i \Phi_j)\ \\
			s& =(j-1)\cdot (N_1-1)+i,\ t=(\tilde{j}-1)\cdot (N_1-1)+\tilde{i}
		\end{align*}
		and the unknown vector ${\bf X}:= [q_t]_{t=1}^{N}$. Then we obtain the linear algebraic system
		\begin{align*}
			\mathcal{A} {\bf X} = {\bf b}.
		\end{align*}
		
\begin{figure}[htbp]
			\centerline{\includegraphics[width=0.7\linewidth]{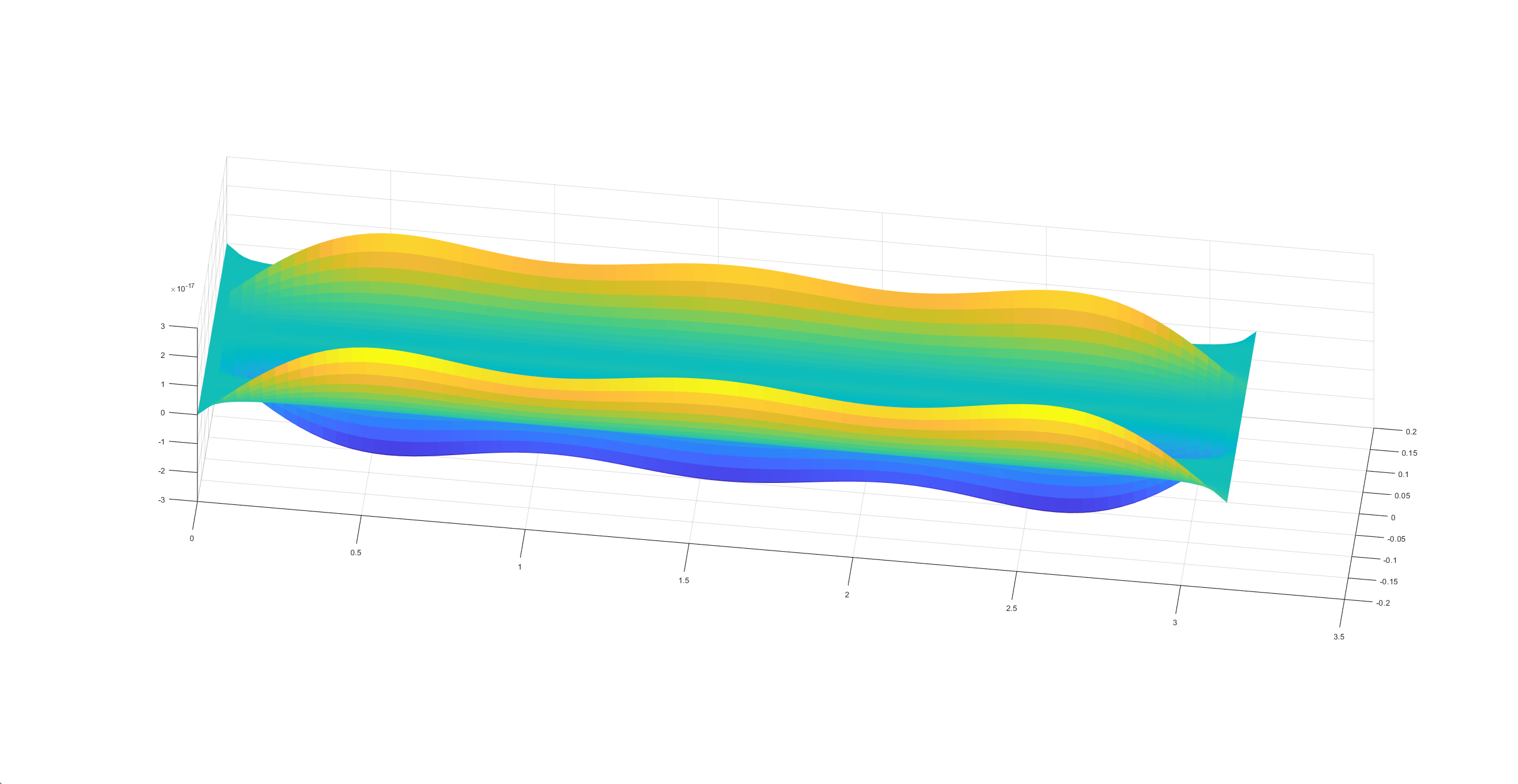}}
			\caption{Solution $U_{-125}$ with multiple zeros at the $x$-direction. Case for $\alpha=-125$ .}\label{sv}
		\end{figure}

		On Fig. \ref{sv} the solution to the linearised problem \eqref{eq:2.1} was performed using the above separation variable method. 
		For implementation of the method and numerical computations the Matlab tool was used.
		
		The nontrivial solution for $U_h$ presented here is for $\alpha=-125$.
		
		\subsection{Triangular mesh case}\label{chapter:numerical}
		In this section the solution to the linearized boundary value problem \eqref{eq:2.1} based on its variational form \eqref{vf4} has been performed.
		The standard finite element method of P4 was applied in order to solve the equation.
		For implementation and computation the FreeFem++ \cite{ff} tool was used.
		
		We set $N=8000$ as the number of iterations and we start with $m=0$ and $\alpha=0$.
		First nontrivial solution for $U_h$ appears at iteration 460 with corresponds to value of $\alpha=-460$.
		The solution has only positive values till iteration 520 ($\alpha=-520$) which means that there is no zeros in the $x$-direction..
		This suggest that the first value of the parameter $m$ is one for $\alpha \in [-520,-460)$

\begin{figure}[htbp]
			\centerline{\includegraphics[width=0.3\linewidth,trim={0 15cm 0  0}]{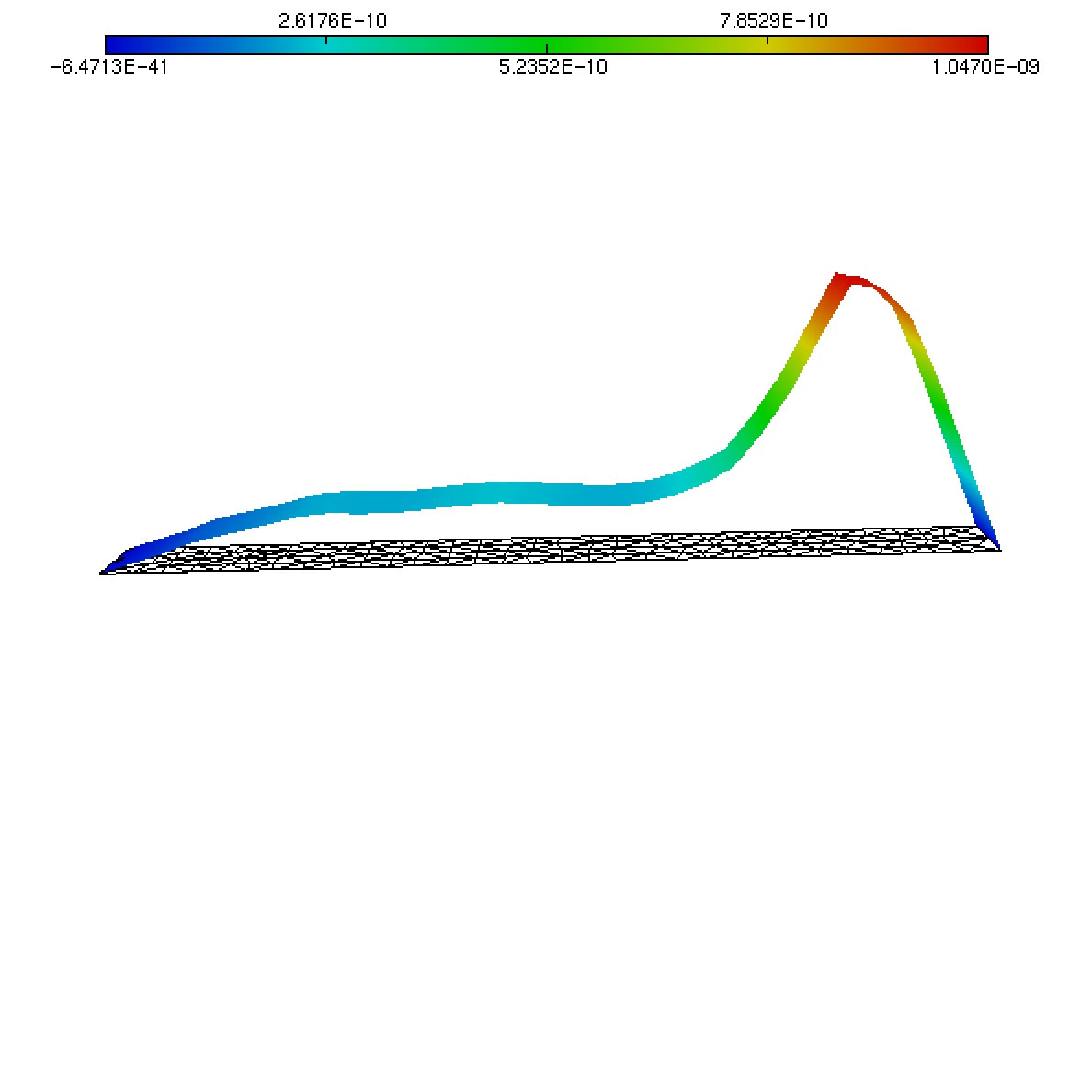}\includegraphics[width=0.3\linewidth,trim={0 15cm 0  0}]{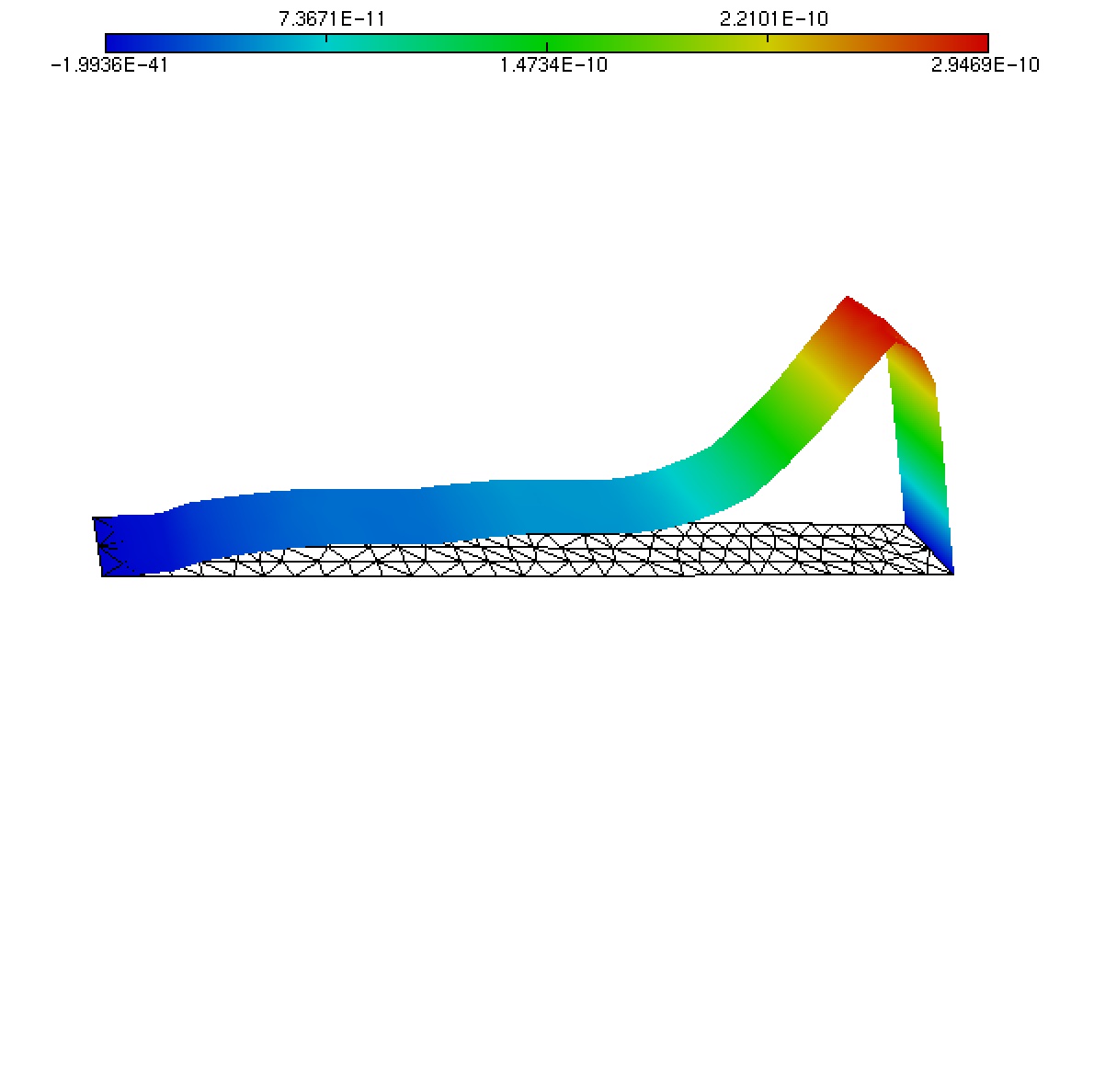}
			}
			\caption{Solution $U_{1,-460}$ with no zeros at the $x$-direction. Case $m=1$ with $\alpha=-520\leq -460=\alpha_1$.}\label{alfa-520-460}
		\end{figure}

		In Fig.  \ref{alfa-520-460}
		some examples of such solutions $U_{1,\alpha}$ where $\alpha=-520\leq -460=\alpha_1$ are presented.
		
\begin{figure}[htbp]
			\centerline{\includegraphics[width=0.3\linewidth,trim={0 15cm 0  0}]{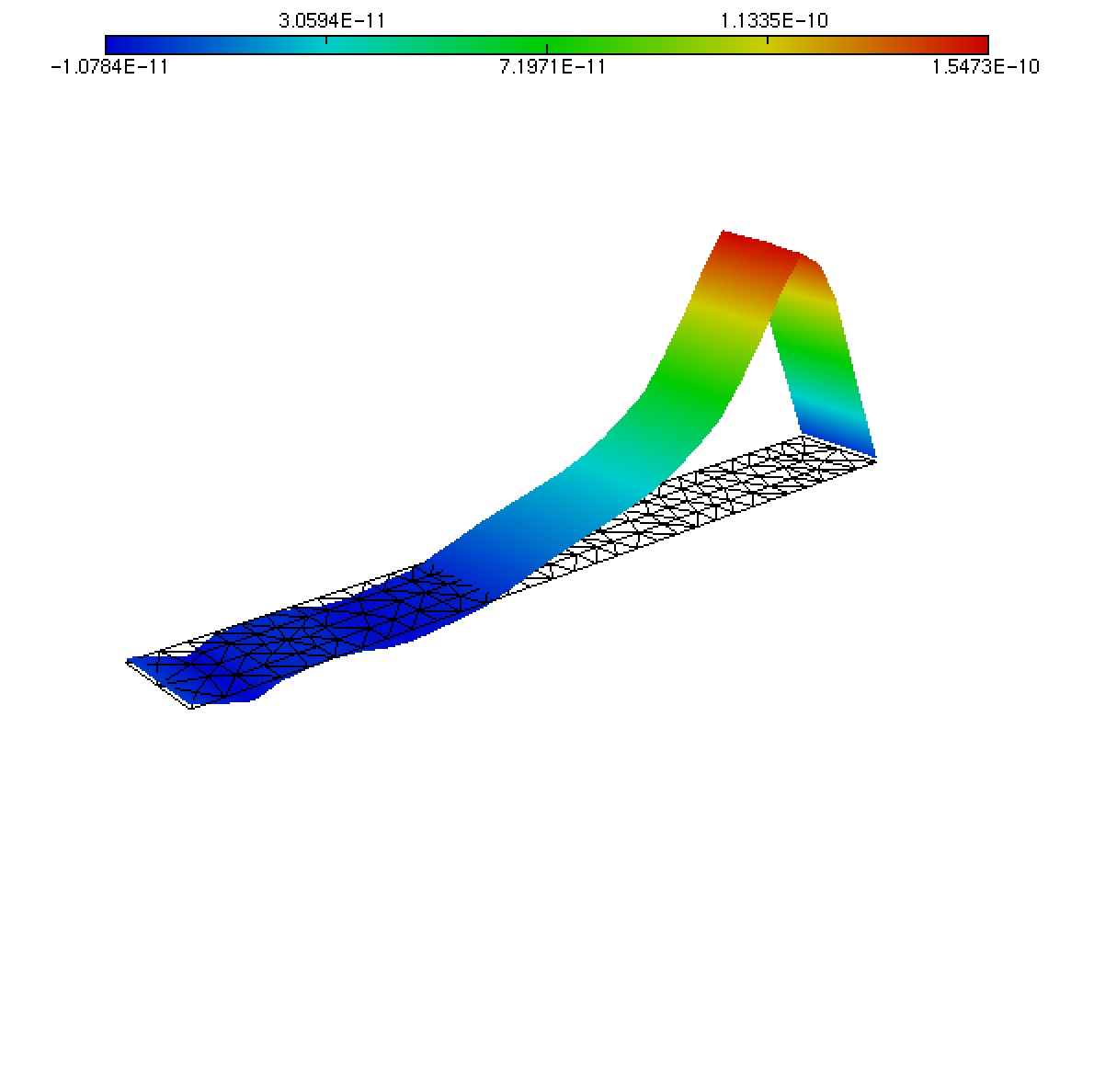}\includegraphics[width=0.3\linewidth,trim={0 15cm 0  0}]{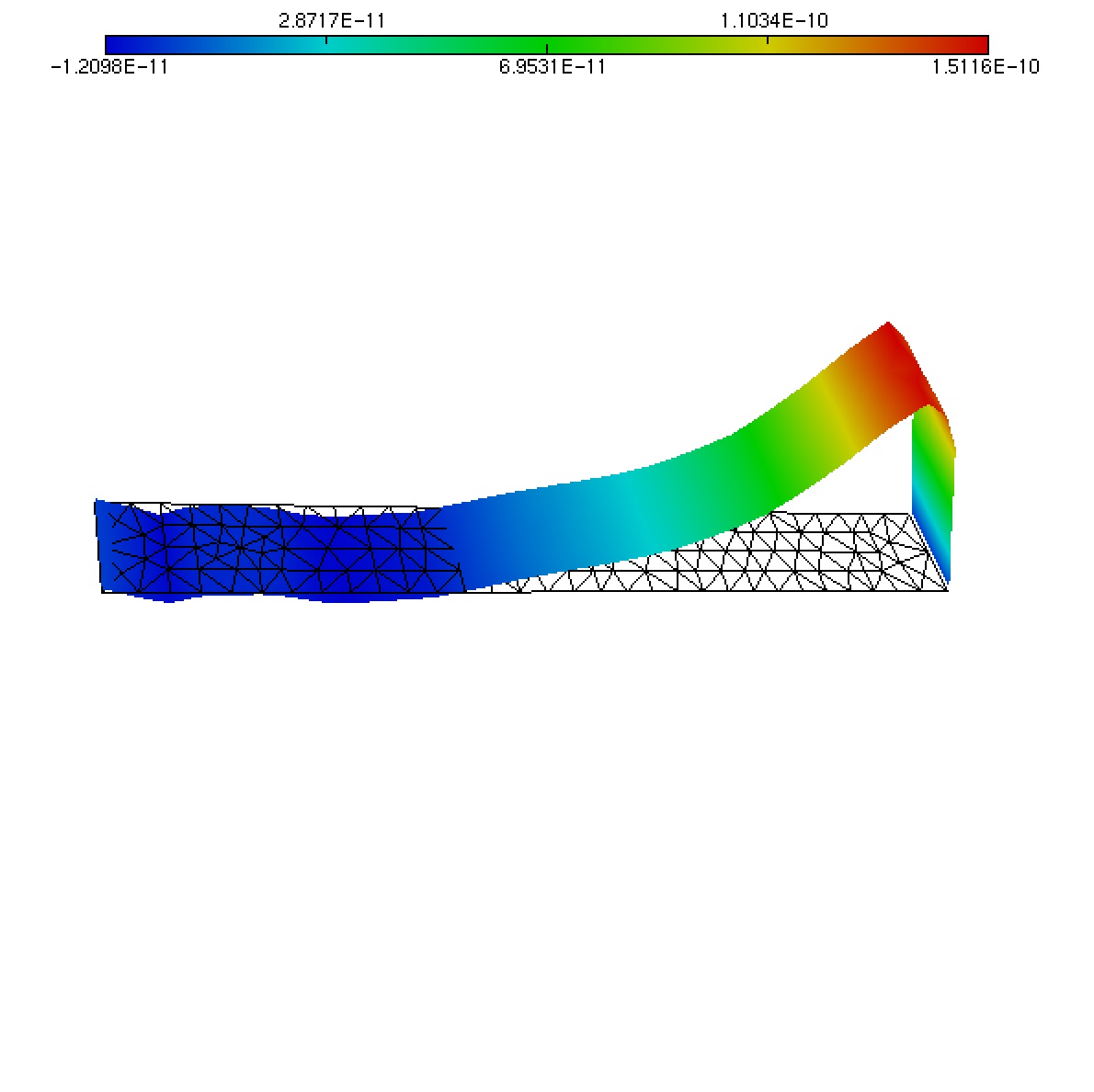}\includegraphics[width=0.3\linewidth,trim={0 15cm 0  0}]{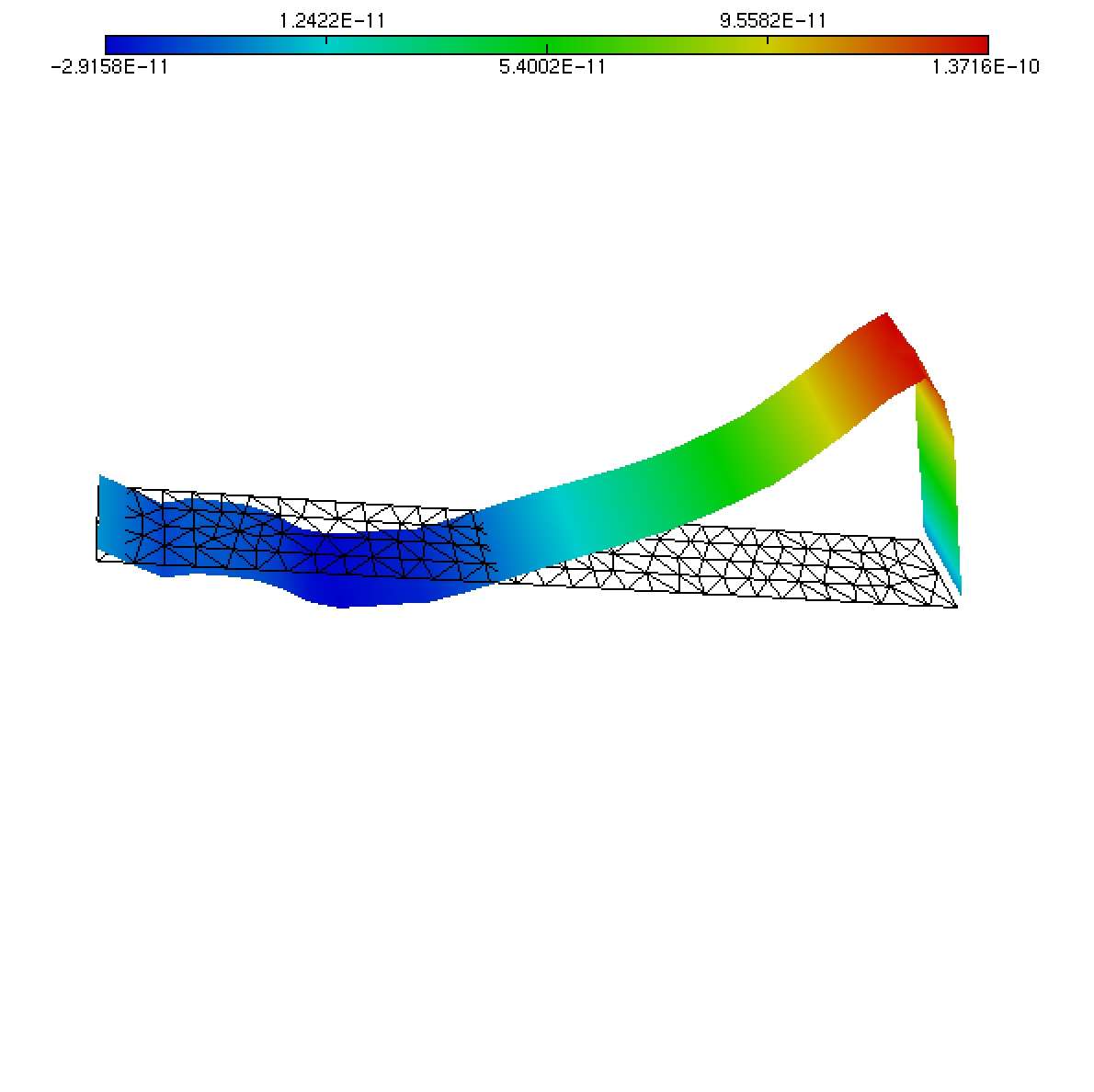}}
			\caption{Solution $U_{2,-520}$ with two zeros at the $x$-direction. Case $m=2$ with $\alpha=-850<-520=\alpha_2$.}\label{alfa-850-520}
		\end{figure}

		In Fig. \ref{alfa-850-520} 
		we can see an example of the solution $U_{m,\alpha}$ with one zero in the $x$-direction. The set of unimodal solutions $U_{m,\alpha}$ with such property was obtained for $\alpha\in [-850, -520]$. In this case we get $m=2$ and $\alpha_2=-520$. 
		
\begin{figure}[htbp]
			\centerline{\includegraphics[width=0.3\linewidth]{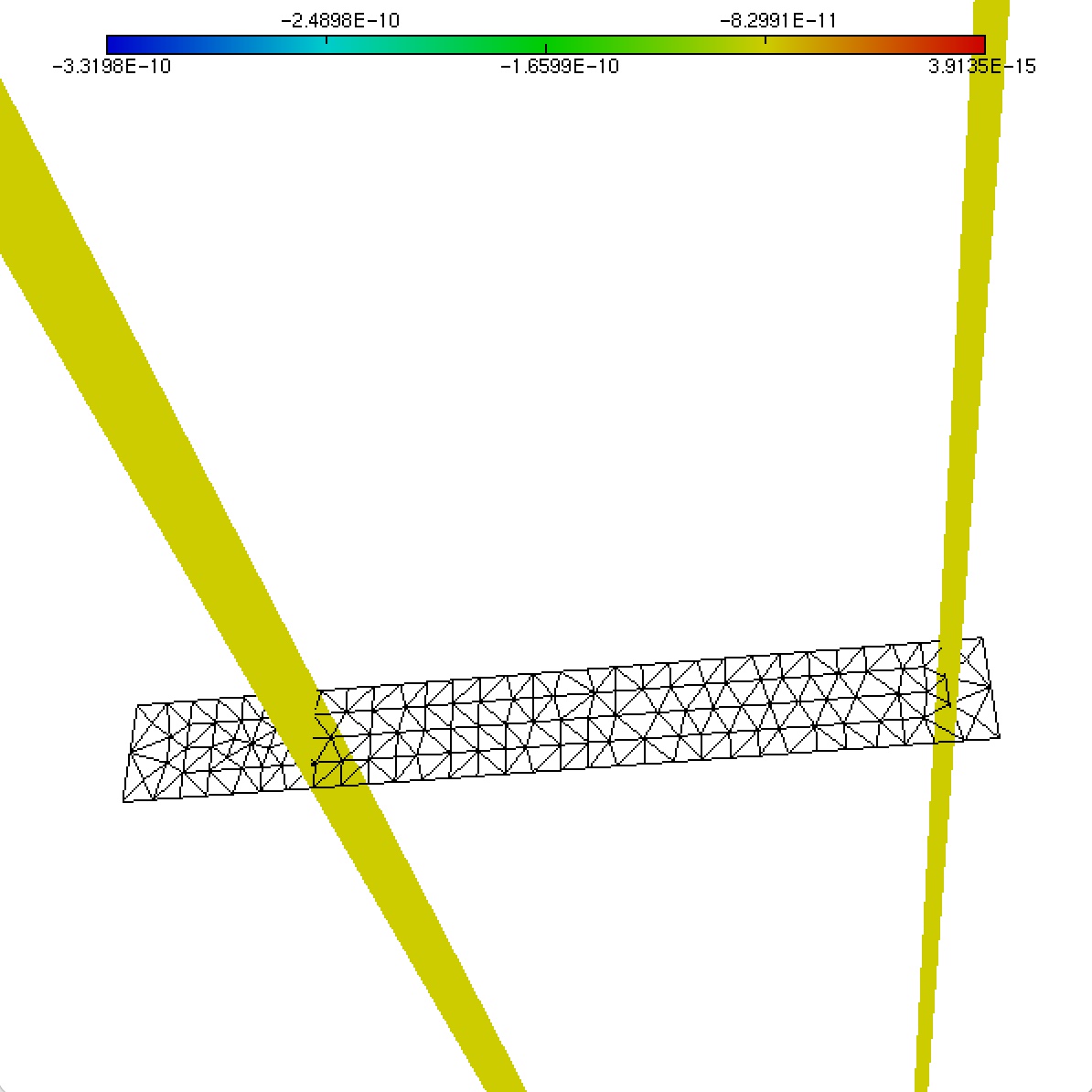}\includegraphics[width=0.3\linewidth]{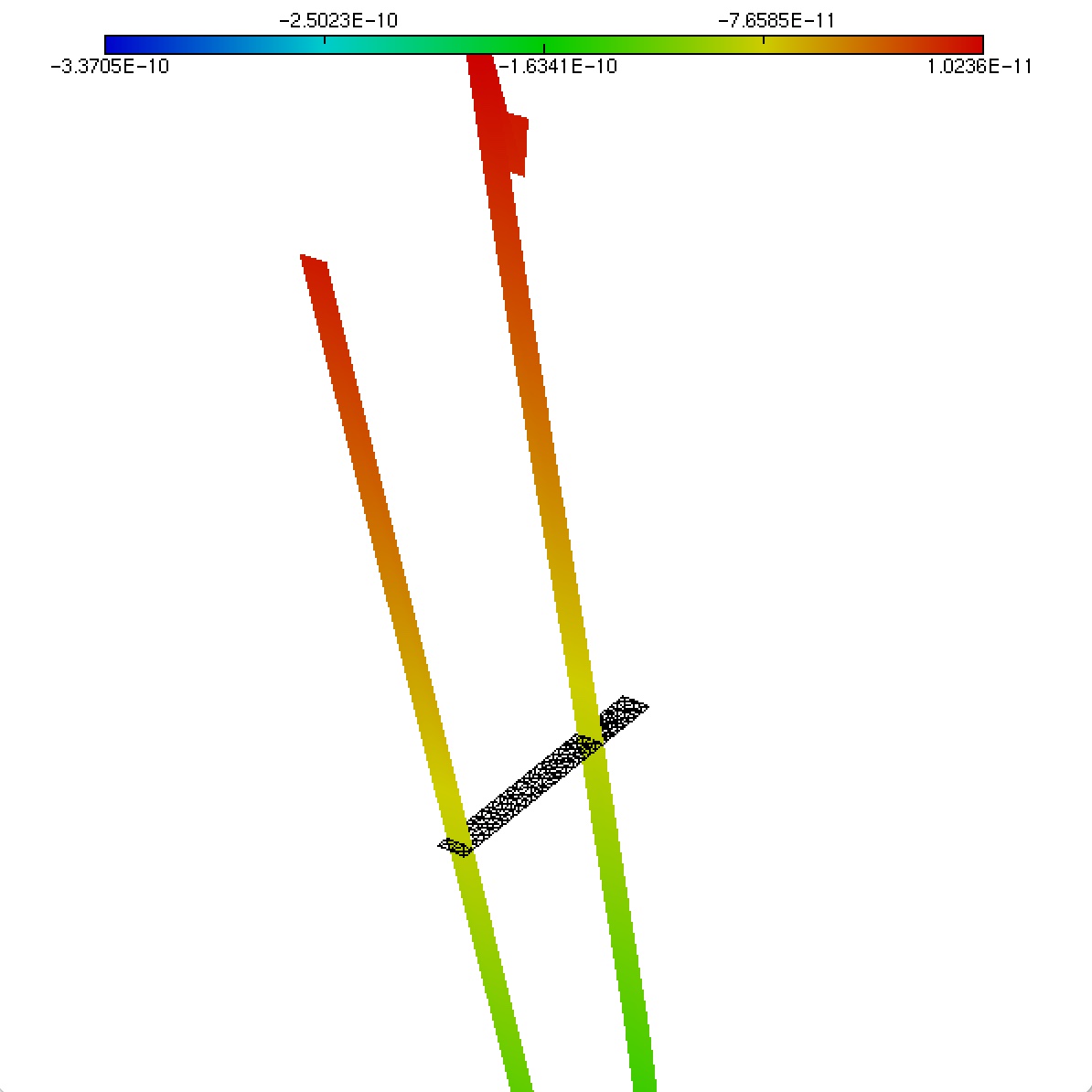}\includegraphics[width=0.3\linewidth]{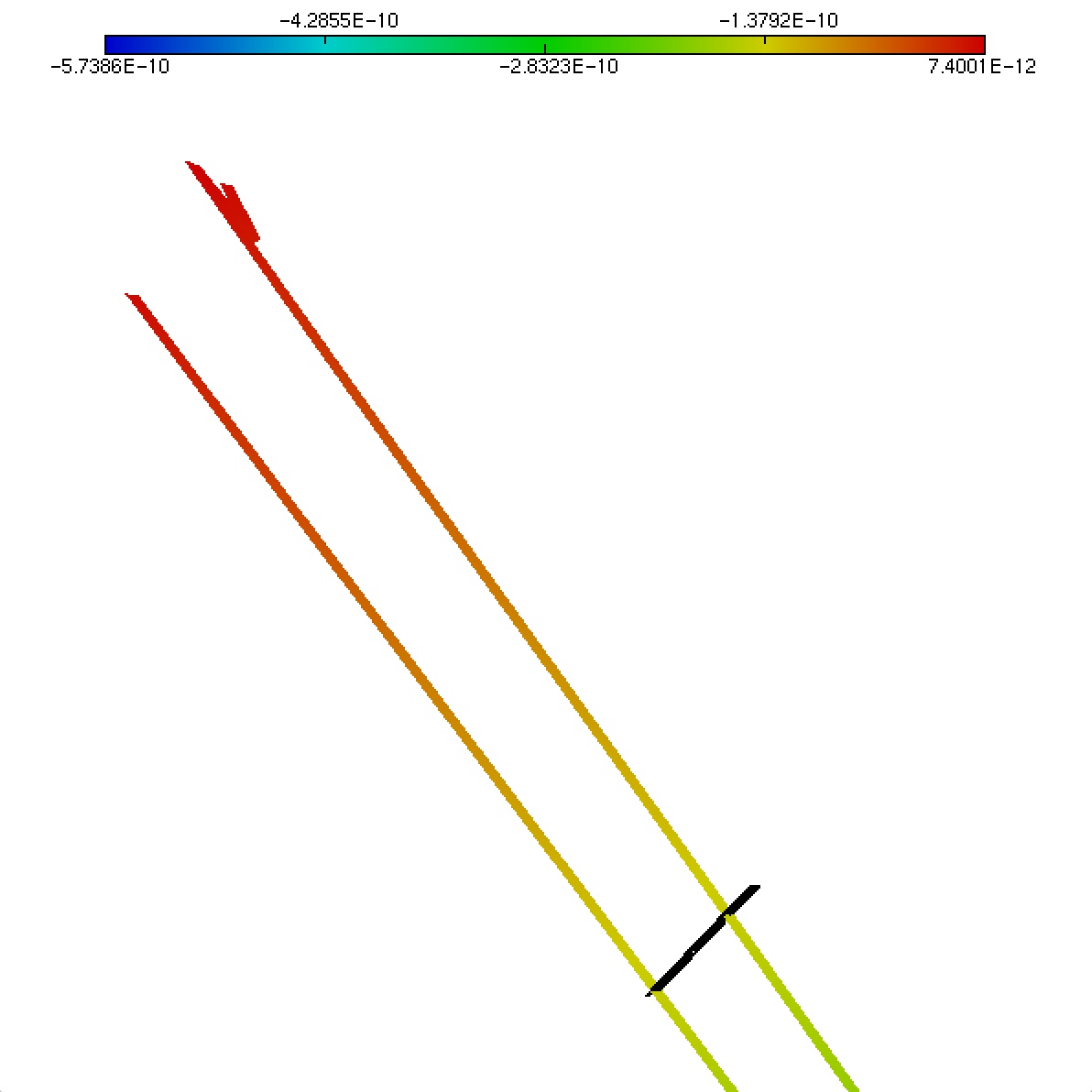}}
			\caption{Solution $U_{3,-2020}$ with four zeros at the $x$-direction. Case $m=3$ and $\alpha_3=-2020$ with $\alpha=-2330<-2020=\alpha_3$.}\label{alfa-2330-2020}
		\end{figure}

		Next example is for $m=3$ and $\alpha_3=-2020$ for which the solution $U_{m,\alpha}$ admits up to 2 zeros in the $x$-direction. This set of solutions is obtained for $\alpha\in[-2330,-2020]$. See 
		Fig.\ref{alfa-2330-2020} 
		for $\alpha=-2020<\alpha_3$ chosen arbitrarily but from given interval.

		\section*{Conclusions}
		In the paper, the numerical solution to the linearized model was considered.
		The main goal was to study behavior of steady state solutions   for different values of the wind velocity given by $\alpha$.
		These are related to asymptotic  [long time] behavior of evolutions driven by \eqref{eq:2.1} as shown in \cite{lasweb}. 
		The reachness of steady states indicates a complex structure of the corresponding attractor  (see \cite{lasweb}) which captures long time behavior of dynamical system defined by \eqref{eq:2.1}.
		In particular, we focus on the study of the linearized system \eqref{eq:8.2}  whose solutions have one to one correspondence to nonlinear solutions of \eqref{eq:8.1}  via the formula \eqref{eq:11}. The corresponding  relation is provided and proved in Lemma \ref{lem22}.
		The observation of the numerical results led us to the conclusion that 
		with the increased values  of $|\alpha|$ the number of solutions in a stationary set increases. 
		The latter explains the effect of the wind on the  flutter  affecting the structure.
		In fact, based on \cite[Corollary 3.6]{lasweb} one expects that  with $|\alpha| \rightarrow \infty$  that number  tends to infinity and this fact was provided by the numerical computations.
		Also P4 finite element method seems to be better approximation for this problem however it needs to be more analyzed.

		\section*{Acknowledgment}
		
		This research was partially founded by the Polish National Science Centre (NCN), grant Opus. Agreement UMO-2023/49/B/ST1/04261.


\begin{thebibliography}{00}
			
			\bibitem{agar1} 
			T.J.A. Agar, 
			``The analysis of aerodynamic flutter of suspension bridges'', 
			Computers and Structures, vol. 30(3), pp. 593--600, 1988.
			
			\bibitem{agar2} T.J.A. Agar,
			``Aerodynamic flutter analysis of suspension bridges by a modal technique'', 
			Engineering structures, vol. 11(2), pp.75--82, 1989.
			

			\bibitem{lasweb}
			D. Bonheure, F. Gazzola, I. Lasiecka, J. Webster,
			``Long-time dynamics of a hinged-free plate driven by a nonconservative force'',
			Ann. Inst. H. Poincar\'e C Anal. Non Lin\'eaire, vol. 39, pp. 457--500, 2022.
			%
			\bibitem{12}
			D. Bonheure, F. Gazzola, E. Moreira dos Santos,
			``Periodic solutions and torsional instability in a nonlinear nonlocal plate equation'',
			SIAM J. Math. Anal., vol. 51(4), pp. 3052--3091, 2019.
			%
			\bibitem{n16}
			I. Chueshov, I. Lasiecka,
			``Von Karman Evolution Equations: Well-posedness and Long Time Dynamics'',
			Springer Science and Business Media, 2010.
			%
			\bibitem{CLW2014} 
			I. Chueshov, I. Lasiecka, J.T. Webster, Attractors for Delayed, Nonrotational von Karman ``Plates with Applications to Flow-Structure Interactions Without any Damping'', Communications in Partial Differential Equations, vol. 39(11), pp. 1965--1997, 2014.

			\bibitem{dowell}
			E.H. Dowell et al., 
			``A Modern Course of Aeroelasticity'', 
			4th edn, Kluwer Academic Publ., Dordrecht, The Netherlands, 2004.

			\bibitem{gaz23}
			F. Gazzola, M. Jleli, B. Samet,
			``A new detailed explanation of the Tacoma collapse and some optimization problems to improve the stability of suspension bridges'',
			Math. Eng., vol. 2(045), 1-35, 2023.
			
			\bibitem{highorder}
			H. Gravenkamp, A. Saputra, S. Duczek. ``High-order shape functions in the scaled boundary finite element method revisited." Archives of Computational Methods in Engineering 28, pp. 473-494, 2021.
			
			\bibitem{ff}
			F. Hecht,
			``New development in FreeFem++'',
			J. Numer. Math., vol. 20(3/4), pp. 251--265, 2012.

			%
			\bibitem{MANSFIELD}
			E.H. Mansfield, 
			``The Bending and Stretching of Plates'',  
			2nd ed., Cambridge University Press, Cambridge, 1989.

			\bibitem{szulc}
			K. Szulc,
			``Numerical solution to the linearized model of a clamped-free plate using nonconforming finite elements'',
			Evolution Equations and Control Theory, 12(6), pp. 1456-1472, 2023.

			\bibitem{21}
			O.C. Zienkiewicz,
			``The Finite Element Method in Engineering Science'',
			Mac Graw-Hill, London, 1971.
			%
			
		\end{thebibliography}
	\end{document}